\documentstyle[11pt]{article}
\input diagrams
\title{The Chow groups and the motive of the Hilbert scheme of points on a surface
}
\author{
Mark Andrea A.  de Cataldo\thanks{
Partially supported by N.S.F. Grant DMS 9701779.}\, 
and Luca Migliorini\thanks{ member of GNSAGA, supported by MURST funds, and by SFB237}
}

\newtheorem{tm}{Theorem}[subsection]
\newtheorem{lm}[tm]{Lemma}
\newtheorem{pr}[tm]{Proposition}
\newtheorem{rmk}[tm]{Remark}
\newtheorem{cor}[tm]{Corollary}
\newtheorem{ex}[tm]{Example}
\newtheorem{defi}[tm]{Definition}

\newtheorem{??}[tm]{Question}

\font\tenmsb=msbm10
\font\sevenmsb=msbm7
\font\fivemsb=msbm5

\newfam\msbfam
\textfont\msbfam=\tenmsb
\scriptfont\msbfam=\sevenmsb
\scriptscriptfont\msbfam=\fivemsb
\def\Bbb#1{{\fam\msbfam #1}}
\font\teneufm=eufm10
\font\seveneufm=eufm7
\font\fiveeufm=eufm5
\newfam\eufmfam
\textfont\eufmfam=\teneufm
\scriptfont\eufmfam=\seveneufm
\scriptscriptfont\eufmfam=\fiveeufm
\def\frak#1{{\fam\eufmfam\relax#1}}

\newcommand\ci{\cite}
\newcommand\s{\sigma}
\newcommand\rat{{\Bbb Q}}

\newcommand\zed{{\Bbb Z}}

\newcommand\pn[1]{{\Bbb P}^{#1}}

\newcommand\blacksquare{{\hspace*{\fill} $\Box$}}

\newcommand\Si{\Sigma}

\newcommand\r{\rho}

\newcommand\J{{\cal J}}

\newcommand\RR{\frak R}

\newcommand\PP{\frak P}

\newcommand\G[1]{\Gamma_{#1}}

\newcommand{\Hi}[1]{X^{[#1]}}
\newcommand{\Sy}[1]{X^{(#1)}}
\newcommand{\Hil}[2]{ X^{[#1]}_{#2}}
\newcommand{\Sym}[2]{X^{(#1)}_{#2}}
\newcommand{\Gn}[2]{\Gamma^{#1}_{#2}}
\newcommand{\prt}{\PP(n)}
\newcommand\sn[1]{X^{#1}}
\newcommand\snn[2]{X^{#1}_{#2}}

\newcommand\Ch[2]{A_{#1}({#2})}
\newcommand\n{\noindent}
\newcommand\bnu{\bf{\nu}}
\newcommand\bmu{\bf{\mu}}
\newcommand\nnn{\succeq  \!\!\!\!\! /}

\begin{document}
\maketitle
\begin{abstract}
We compute  the  Chow  motive  and the Chow groups 
with rational coefficients of the Hilbert scheme
of points on a smooth algebraic surface.

\end{abstract}

%\tableofcontents

\parindent0pt

\section{Introduction}
\label{intro}
In this paper we compute  the Chow  motive and the Chow groups 
with rational coefficients of the Hilbert scheme
$ X^{ [n] } $
of $n$ points on a nonsingular algebraic surface $X$.
The result holds over an arbitrary field.

The space $ X^{ [n] } $
is a remarkably rich object which finds itself, for reasons that are perhaps still
mysterious, at the croassroads of geometry, representation theory and mathematical physics.

In the context of  testing the $S$-duality conjecture, using G\"ottsche's
calculation of the Betti numbers of $ X^{ [n] } $, Vafa and Witten
suggested, for reasons stemming from orbifold cohomology, that
the singular cohomology groups of these Hilbert schemes should be naturally
linked to infinite dimensional graded Lie algebras. This fact was 
firmly established by Nakajima and Grojnowski, independently. 
Some of the contributors to this circle of ideas are: Fogarty \ci{fo}, Briancon \ci{br},
Iarrobino \ci{ia}, 
Ellingsrud-Stromme \ci{e-s1}, G\"ottsche \ci{go}, G\"ottsche-Sorgel
\ci{go-so}, Cheah \ci{ch}, Vafa-Witten \ci{va-wi}, Nakajima \ci{nakp}, Grojnowski
\ci{gr}.
% Lehn \ci{le}. 
The reader is referred to the beautiful lectures \ci{nak}
and to our paper \ci{decatmig} for background, further results and references to the literature.

In our paper \ci{decatmig} we proved, directly,  a precise form  of 
the decomposition theorem for the so-called Hilbert-Chow map.
We detected, via the action of the Lie algebra, 
certain   
subvarieties of $ X^{ [n] } $  carrying 
 the topological information necessary to understand the additive structure of singular cohomology.
Subsequently, we  looked directly at certain related correspondences.
 
In this paper we show how these correspondences identify the Chow
groups
   with rational coefficients
 of $ X^{ [n] } $ with  the Chow groups with rational coefficients 
 of a certain collection of 
 products of symmetric products of the surface $X$. 
 Using this result we  determine the Chow motive with rational coefficients
 of $ X^{ [n] } $.

Voevodsky's theories of motivic cohomology  incorporate the formalisms of Tate twists and shifts of complexes.
Once the Chow motive  of $ X^{ [n] } $ is computed, it is easy to
compute the motive in these more general theories. Our result, read in these theories,
 harmonizes the shifts  present in the  decomposition theorem mentioned above and the 
 Tate twists in the mixed Hodge structure of $ X^{ [n] } $.
There is also a generating function for the motives in question.

 The main tool used in our calculation is the Gysin formalism of 
 Fulton-MacPherson. None of the tools used in the various derivations
 of the structure of the  singular cohomology  of $ X^{ [n] } $ 
 in the literature seems applicable in the context of cycles. In fact, this cycle-theoretic approach can be easily supplemented to give another proof of G\"ottsche's formula for the Betti numbers.
  The fact, due to Ellingsrud and Stromme, that punctual Hilbert schemes admit affine cellular decompositions is essential to our approach.

 The outline of the paper is as follows.
 $\S$\ref{cubo} is devoted to fixing the notation and to introducing 
 the varieties naturally associated
 with  $X^{[n]}$  and  its natural stratification given by  partitions.
 $\S$\ref{intersection}  reviews  the basic results from intersection
 theory that we need.  $\S$\ref{correspondences}  reviews  well-known facts
 concerning the Gysin formalism of correspondences in a ``non-complete" situation.
  $\S$\ref{cvrg}
 discusses correspondences in this 
 ``refined" formalism. 
 $\S$\ref{coc} discusses the composition of correspondences in the context
 of the refined formalism.  $\S$\ref{qv} overviews the situation
 for quotient varieties. 
 In $\S$\ref{mgg} we define the
 natural map $\hat\Gamma_*$ via the correspondences $\hat\Gamma$. 
 The main result of the paper is Theorem \ref{maintm}, which states that
 $\hat\Gamma_*$ is an isomorphism. The injectivity statement is  Corollary 
 \ref{injj}.
  The surjectivity statement, Corollary \ref{surjquot}, is the heart of the present paper.
  Corollary \ref{groth} and Remark \ref{equi} are  standard consequences for the Grothendieck
 groups also in the equivariant context.
 $\S$\ref{motive} is devoted to the identification, Theorem \ref{motiff}, of the motive of the Hilbert scheme
 $X^{[n]}$ with a collection of  motives of products of symmetric products of the surface $X$. 
 Theorem \ref{vmot} is the translation of Theorem \ref{motiff} into
 Voevodsky's categories. We give a ``generating function" for this motivic structure
 at the end of $\S$\ref{projd}.

 \bigskip

 {\bf Acknowledgments.} Some ideas of N. Fakhruddin have
 influenced our approach to this question. 
  We have asked him to co-author the present paper, but
 he has declined.
 The first author would like to thank 
 Hong Kong University for its very generous hospitality while 
 portions of this work
 have been carried out. We wish to thank M. Kumar, J.I. Manin and A. Vistoli for useful conversations.
The commutative diagrams are produced by Paul Taylor's \texttt{diagrams}
macro package, available in CTAN in
\texttt{TeX/macros/generic/diagrams/taylor}. 
 
\section {Some special varieties and fibrations}
\label{definitions}
In these section we review a few definitions and introduce the correspondences which 
will play a major role in the paper. Proofs or references for all the results stated here can be found in our previous paper \cite{decatmig}.
Let $X$ be an irreducible  quasiprojective nonsingular surface 
defined over an algebraically closed field, $\Sy{n}$ be 
its n-th symmetric product, $\Hi{n}$  be the Hilbert scheme of 0-dimensional subschemes of 
$X$ of length $n$, $\pi: \Hi{n} \to \Sy{n}$ be  the Hilbert-Chow morphism. 
It is well known that $\Hi{n}$ is nonsingular. $\prt$ denotes the set of partitions of $n$ and $p(n)$ its cardinality. 
If $\nu \in \prt$, we denote by $l(\nu)$ its length, and define $\Sym{n}{\nu}$ 
to be the locally closed subset of points in $\Sy{n}$ of the type 
$\nu_1x_1+ \cdots  + \nu_{l(\nu)}x_{l(\nu)}$,
 with $x_h \in X$ and $x_i \neq x_j$ for every $i\neq j$. Similarly, 
$\Hil{n}{\nu}: = ( \pi^{-1}(\Sym{n}{\nu}) )_{red}$. Let $\overline \Sym{n}{\nu}$ 
denote the closure of the stratum. It can  be proved that 
$\overline \Hil{n}{\nu}= ( \pi^{-1}(\overline \Sym{n}{\nu}) )_{red}$.
If $\nu=1^{a_1}\cdots n^{a_n}$, then   the finite group 
$\Si_{\nu}:=\Si_{a_1}\times \cdots \times \Si_{a_n}$ 
acts naturally on $\sn{l(\nu)}$. 
The quotient $\Sy{\nu}$ is isomorphic to 
$\Sy{a_1}\times\cdots\times\Sy{a_n}$.
For the sake of uniformity of notation we shall denote
$\sn{l(\nu)}$ by $\sn{\nu}$.
There is a natural  $\Sigma_{\nu}$-invariant map $\nu:\sn{\nu}\to \Sy{n}$ 
whose image is $\overline \Sym{n}{\nu}$. This map 
  descends  to a map 
which we denote by the same symbol $\nu: \Sy{\nu} \to \Sy{n}$.
This map is the normalization map of $\overline X^{(n)}_{\nu} $.
We define correspondences $\Gn{\nu}{}$ and ${\hat{\Gamma}}^{\nu}$ 
as follows:
$$\Gn{\nu}{}=\{\, (x_1,\cdots,x_{l(\nu)}, \J)\in \sn{\nu} \times \Hi{n}:\pi(\J)=
\nu_1x_1+ \cdots + \nu_{l(\nu)}x_{l(\nu)} \, \}\simeq (\sn{\nu} \times_{\Sy{n}} \Hi{n})_{red}. $$
Denote by $\sn{\nu}_{reg}$  (resp.  $\sn{(\nu)}_{reg}$) the set of points of $\sn{\nu}$ (resp. $\sn{(\nu)}$) strictly of type $\nu$.
\begin{rmk}
{\rm The restriction $\Gamma^{\nu}_{reg}$ of $\Gn{\nu}{}$ to $\sn{\nu}_{reg}$ is a Zariski locally 
trivial fibration with irreducible fibers
and it is open and  dense in $\Gn{\nu}{}$ (see \cite{decatmig} Lemma 3.6.1 and
Remark \ref{fibre}). 
It follows that   $\Gn{\nu}{}$ is irreducible of dimension $n+l(\nu)$.}
\end{rmk}

The projection $p_1$ will be denoted by $p_{\nu}$
in order to emphasize the dependence on $\nu$. 
The correspondence $\Gn{\nu}{}$ is  invariant under the action of $\Si_{\nu}$ on the first factor of the product. We can  therefore define $\hat \Gn{\nu}{} : =\Gn{\nu}{}/\Si_{\nu}$.
We summarize these constructions with the diagram:

%\begin{diagram}
%\Gn{\nu}{}            &               &                 &                   &  \\
%\dTo^{p_{\nu}}         & \rdTo^{q'}\rdTo(4,2)^{P}&     &                   & \\
%\sn{\nu}        &               & \hat\Gn{\nu}{}     & \rTo^{P}& \Hi{n} \\
%                   & \rdTo^q       & \dTo^{\hat p_{\nu}} &   & \dTo^{\pi} \\
%                   &            & \Sy{\nu} & \rTo^{\nu}& \Sy{n} \\
%\end{diagram}
\begin{diagram}
\Gn{\nu}{}       &\rTo^{q'}              & \hat\Gn{\nu}{}     & \rTo^{P}& \Hi{n} \\
\dTo^{p_{\nu}}       &        & \dTo^{\hat p_{\nu}} &   & \dTo^{\pi} \\
\sn{\nu}            &  \rTo^q          & \Sy{\nu} & \rTo^{\nu}& \Sy{n} \\
\end{diagram}
where $q$ and $q'$ are  the quotient maps by the action of $\Si_{\nu}$.
The composition $P\circ q'$ will still be denoted by $P$.

%{\bf QUESTO E' NUOVO, RENDE PIU' SNELLA LA SCRITTURA DI ALCUNE DIMOSTRAZIONI}
\begin{rmk}
\label {identify}
{\rm The restriction of $P$ to $\hat \Gn{\nu}{reg}=\Gn{\nu}{reg}/\Si_{\nu}$ identifies  $\hat \Gn{\nu}{reg}$ with the locally closed stratum
$X^{[n]}_{\nu}$}.
\end{rmk}
We  introduce ${\cal X}= \coprod_{\nu \in \prt}\sn{\nu}$ and
 $\hat{\cal X}= \coprod_{\nu \in \prt}\Sy{\nu}$,
 $\Gamma= \coprod_{\nu \in \prt}\Gn{\nu}{}$,  
$\hat \Gamma= \coprod_{\nu \in \prt}\hat \Gn{\nu}{}$
and define $p:\Gamma \to {\cal X}$ (resp. $\hat p: \hat \Gamma \to \hat {\cal X}$)
to be the map induced by  all the maps $p_{\nu}$ (resp. $\hat{p}_{\nu}$).
 
The set of partitions has a natural partial order, which reflects the incidence relations of the strata:
\begin{defi}
Let $\nu,\,\mu \in \prt$. We say that $\nu \succeq \mu$ 
if there exists a decomposition $\{I_1,\cdots,I_{l(\mu)} \}$ of the set $ \{1,\cdots l(\nu)\}$
such that $ \mu_1= \sum_{i \in I_1} \nu_i, \cdots , \mu_{l(\mu)}= \sum_{i \in I_{l(\mu)}} \nu_i $.
\end{defi}
It is easily seen that $\mu \succeq \nu$  if and only if
$\Sym{n}{\nu} \subseteq \overline \Sym{n}{\mu}$.
Fix a total order $ \geq$ on $\prt$ which is compatible with
$\succeq$. For any $\mu \in \prt$ we have the open subsets 
$\Sym{n}{\geq \mu}:= \coprod_{\nu \geq \mu}\Sym{n}{\nu} \subseteq X^{(n)}$.   
Similarly  for $\Sym{n}{>\mu}$. Correspondingly, we have  
open sets ${\cal X}_{\geq \mu}$,  ${\cal X}_{> \mu}$, $\G{\geq \mu}$, $\G{> \mu}$, 
$\Hil {n}{\geq \mu}$ and $\Hil{n}{ > \mu}$
obtained by base change followed by reduction. We  have the corresponding quotients
$\hat{{\cal X}}_{\geq \mu}$,  $\hat{{\cal X}}_{> \mu}$, 
$\hat{\Gamma}_{\geq \mu}$, $\hat{\Gamma}_{ > \mu }$.
Similarly, with  the symbol $\geq$ replaced by $\succeq$.

The imbedding $\Sym{n}{\mu}\to \Sym{n}{\geq \mu}$ is closed. We have 
corresponding closed imbeddings
${\cal X}_{\mu} \to {\cal X}_{\geq \mu}$ and $\G{\mu} \to \G{\geq \mu}$ 
obtained by base change followed by reduction. 
Note that $\G{\mu}\neq \Gn{\mu}{}$.
We thus have a diagram:
\begin{diagram}
\label{cubo}
\G{\mu}              &                  &  \rTo^{P}             &                & \Hil{n}{\mu}      &                     &                    \\
                     &     \rdTo        &                       &                & \vLine ^{\pi}     & \rdTo               &                    \\
\dTo^{p_{\geq \mu}} &                  & \G{\geq \mu}          & \rTo^P         & \HonV             &                     & \Hil{n}{\geq \mu}  \\ 
                     &                  &     \dTo^{p_{\geq \mu}} &                & \dTo              &                     &                    \\
 {\cal X}_{\mu}      &   \hLine         &  \VonH                & \rTo           & \Sym{n}{\mu}      &                     & \dTo_{\pi}         \\
                     &  \rdTo           &                       &                &                   &       \rdTo         &                     \\  
                     &                  & {\cal X}_{\geq \mu}   &                & \rTo              &                     &   \Sym{n}{\geq \mu} . \\
\end{diagram}

\section{Review of Intersection Theory}
\label{intersection}
\bigskip
Our unique reference for this section is \cite{fulton}.
In what follows, the Chow groups $\Ch{\ast}{X}$ of an algebraic scheme $X$ over a field
are always 
taken with rational coefficients, 
even though many results hold with integer coefficients.
We recall that given a regular imbedding $i:X \to Y$ of codimension $d$
with normal bundle $N_XY$ and a morphism $f:Y' \to Y$ from a pure $l$-dimensional variety
$Y'$, there are refined Gysin homomorphisms 
$i^!:\Ch{\ast}{Y'}\to \Ch{\ast-d}{X'}$ where 
$X'=X \times_Y Y'$ . The construction of $i^!([Y'])$ goes as follows: 
the map $X' \to Y'$ is a closed imbedding and the normal cone $C_{X'}Y'$ is a pure $l$-dimensional 
subscheme of the pullback  of $N_XY$ to $X'$ ;
its cycle class is therefore equivalent to the flat pullback of a unique cycle class  
$i^!([Y'])\in \Ch{l-d}{X'}$. 
%Using the notion of Segre class \ci{fulton}, $\S4$,
% $ i^!([Y'])$ is the $l-d$-dimensional part of $c(N_X Y_{|Y'})\cap s(X',Y')$.
To define  $i^!(\alpha)$ for $\alpha \in \Ch{\ast}{Y'}$ one replaces 
$Y'$ with the support of $\alpha$ and maps the resulting class to $\Ch{\ast}{Y'}$.

\medskip

Let $f:X\to Y$ be a morphism from a scheme $X$ to a nonsingular variety $Y$. 
The graph morphism $\gamma_f:X\to X\times Y$ is a regular imbedding. Given 
maps $X'\to X$ and $Y'\to Y$ there is a refined Gysin morphism 
$\gamma_f^!:\Ch{k}{Y'}\otimes \Ch{l}{X'}\to \Ch{k+l-dimY}{X'\times_Y Y'}$.
%If $\alpha\in \Ch{\ast}{X'}$ and $\beta\in \Ch{\ast}{Y'}$
%the class $\gamma_f^!(\alpha\otimes\beta)$ is also denoted $\alpha \times_f\beta$.
If $X'=X$, $Y'=Y$ and both maps are the identity map, then we will denote
$\gamma_f^!(\alpha \otimes \beta)$,
 the image
of $\alpha \otimes \beta$
 via  the morphism
$\Ch{k}{Y} \otimes \Ch{l}{X} \to \Ch{k+l-dimY}{X}$,
with the more suggestive piece of notation $f^*(\alpha)\cap\beta$.
If $X'= X$, then we will use the notation
$f^!(\alpha)\cap\beta$ for the refined  intersection.

\begin{rmk}
\label{propint}
{\rm Let $\alpha \in Z_k(Y'),\, \beta \in Z_l(X')$ be irreducible cycles. If the irreducible components
$\xi_i$ of $|\alpha|\times_Y |\beta|$ have the expected dimension
$k+l-dimY$, then $\gamma_f^!(\alpha \otimes \beta)$ is a linear combination
with strictly positive coefficients of the cycles $\xi_i$; 
see \ci{fulton}, $\S7.1$.
In particular, if $Y'$ is a closed subvariety of $Y$ and 
$f^{-1}(Y')$ is irreducible of dimension $dimY'+dimX-dimY$, 
then  $f^!([Y'])\cap [X]$
is a positive multiple of $[f^{-1}(Y')]$.}
\end{rmk}
\medskip

%The following easy lemmata, whose proof is left to the reader,
%  are used in the proof of Lemma
%\ref{fibration}.

%\begin{lm}
%\label{commexseq}
%Let

%\begin{diagram}
%A       & \rTo^p   &  B        & \rTo^q         & C          & \rTo    &  0  \\
%\uTo^f  &          & \uTo^g    &                 & \uTo^h    &         &     \\
%D       & \rTo^r   &  E        & \rTo^s          & F         & \rTo    &  0 \\
%\end{diagram}

%be a commutative diagram of exact sequences of abelian groups. 
%If $f$ and $h$ are surjective, then so is 
%$g$.
%\end{lm}

%\begin{lm}
%\label{oneall}
%Let $X$ be an irreducible  nonsingular variety, 
%$p: \Gamma \to X$ be a morphism which is Zariski locally trivial  
%with fiber $F$. Assume that $F$  is an algebraic scheme for which rational and %algebraic equivalence coincide.
%Let $\alpha \in A_{\ast}(\Gamma)$ be such that $\alpha_x \sim_{rat} 0$ for some %point $x \in X$
%(here $\alpha_x$ is the specialization of $\alpha$ at the fiber $\Gamma_x$ %(\ci{fulton}, $\S10.3$) and is independent
%of the local trivialization of the morphism around $x\in X$).
%Then $\alpha_{x'}\sim_{rat} 0$ for every $x' \in X$.
%\end{lm}

%\medskip

The following  well-known lemma will be used in the sequel of the paper.
We omit the simple proof which  follows easily from the case of a trivial fibration (see \ci{fulton}, Ex. 1.10.2) and by noetherian induction
using the basic properties of the refined Gysin formalism.
\begin{lm}
\label{fibration}
Let $X$ be a nonsingular irreducible variety and $p:\Gamma \to X$ a Zariski locally 
trivial fibration with fiber $F$ admitting a cellular decomposition.
Suppose $\{ \alpha_i \}_{i \in I}$ is a set of classes in
$\Ch{\ast}{\Gamma}$ whose restrictions generate $\Ch{\ast}{p^{-1}(x)}$ for every $x \in X$.
Then  $\{ \alpha_i \}_{i \in I}$ is a set of generators of the
$\Ch{\ast}{X}$-module  $\Ch{\ast}{\Gamma}$. In other words, 
the map $\Phi: \Ch{\ast}{X}^{\oplus I} \to \Ch{\ast}{\Gamma}$, defined by
$\Phi(\{ \beta_i\})= \sum_{i}{ p^*(\beta_i)\cap \alpha_i}$, is surjective.
\end{lm}

\section{Correspondences}
\label {correspondences}

\subsection{Correspondences via refined Gysin maps}
\label{cvrg}
The standard reference is \ci{fulton}, Remark 16.1,  $\S 6$ and  $\S 8$. Consider a diagram
\begin {diagram}
\Gamma & \rTo^{p_2} & X_2 \\
\dTo^{p_1}&  & \\
X_1  &  &   \\
\end{diagram}
with $X_1$ nonsingular and $p_2$ proper.
In section \ref{intersection} we define  a map
 \begin{diagram}
\Ch{k}{X_1} \otimes \Ch{l}{\Gamma} &
                          \rTo^{p_1^*( - )\cap( - )}& 
                                                 \Ch{k+l-dimX_1}{\Gamma}&
                                                                          \rTo^{p_{2*}} & \Ch{k+l-dimX_1}{X_2}.\\
\end{diagram}

The case we will mostly use is 
$$
p_{2*}(p_1^*( - )\cap [\Gamma]):
\Ch{k}{X_1} \longrightarrow \Ch{k+dim \Gamma-dimX_1}{X_2}, 
$$
where $[\Gamma]$ is some  fixed cycle. 
%Tipically, $[\Gamma]$
%will be the fundamental cycle associated with $\Gamma_{red}$.

Most important for us will be the case of relative correspondences 
of the type
\begin {diagram}
\Gamma & \rTo^{p_2} & X_2 \\
\dTo^{p_1}& \Box & \dTo \\
X_1  & \rTo^f & Y,    \\
\end{diagram}
where $X_1$ is a nosingular variety, $f$ is a proper morphism
and $\Gamma=X_1\times_Y X_2$.
If $V \to Y$ is a morphism, then  we have a similar diagram by pulling back to $V$.
%\begin {diagram}
%\Gamma_V & \rTo^{p_2} & X_{2,V} \\
%\dTo^{p_1}&  & \dTo \\
%X_{1,V}  & \rTo & V   \\
%\end{diagram}
Set, for brevity, $X_{1,V}=X_1 \times_Y V$ and similarly  $X_{2,V}=X_2 \times_Y V$
and $\G{V}=\Gamma \times_Y V$. Note that, 
by the transitivity property of the fibre product,
$\G{V}=X_{1,V} \times_V X_{2,V}$.
Fix a cycle $[\Gamma]$.  The refined intersection product gives a map
%$\gamma_{p_1}^!( - ???? \otimes [\Gamma]):$
$p_1^!( - )\cap [\Gamma] :\Ch{\ast}{X_{1,V}}\to \Ch{\ast}{\G{V}} $.

\begin{rmk}
{\rm It is important that we do not assume that $X_1$ is complete, i.e.
that we do not assume 
that $p_2:X_1 \times X_2 \to X_2$ is proper and that we do 
not factorize through
$\Ch{\ast}{X_1 \times X_2}$. Otherwise,   we would loose too much  information: 
e.g. the diagonal in ${\Bbb A}^1 \times {\Bbb A}^1$ is zero in $A_1({\Bbb A}^1 \times {\Bbb A}^1)$, 
but it induces the identity map via the formalism of refined intersections.}
\end{rmk}

Let  $i: Z \to Y$ be a closed imbedding and 
$j:U:=Y \setminus Z \to Y$ be  the resulting open imbedding.
\begin{lm}
\label{exact}
Let $[\Gamma]$ be any cycle in $A_{\ast}(\Gamma)$.
The diagram 
\begin{diagram}
\Ch{\ast}{X_{1,Z}} & \rTo^{p_1^!( - )\cap [\Gamma]} & \Ch{\ast}{\Gamma_Z} & \rTo^{p_{2*}}& \Ch{\ast}{X_{2,Z}} \\
\dTo^{i_*}  &              & \dTo                                       &                                      & \dTo  \\
\Ch{\ast}{X_1}            & \rTo^{p_1^*( - )\cap [\Gamma]}           & \Ch{\ast}{\Gamma}            & \rTo^{p_{2*}}& \Ch{\ast}{X_{2}} \\
\dTo^{j^*}  &              & \dTo                 &              & \dTo  \\
\Ch{\ast}{X_{1,U}} & \rTo^{p_1^*( - )\cap [\Gamma_U]} & \Ch{\ast}{\Gamma_U} 
                                                                                & \rTo^{p_{2*}}& \Ch{\ast}{X_{2,U}} \\ 
\dTo         &             &   \dTo               &              & \dTo  \\
0             &             &    0                 &              &  0      \\              
\end{diagram}
is commutative and the columns are exact. 
\end{lm}
{\em Proof.} The exactness of the columns stems from
\ci{fulton}, Proposition 1.8. The commutativity of the first column of diagrams
stems from \ci{fulton}, Proposition 8.1.1.(c)  and Theorem 6.2. The commutativity
of the second column is obvious. 
\blacksquare

\bigskip
\begin{rmk}
\label{cmn}
{\rm Note that, even if $\Gamma$ is  reduced,  $\Gamma_Z$ may fail to be so. 
However, this causes no trouble in view of the obvious canonical isomorphism
$A_{\ast}(-_{red}) \to A_{\ast}(-)$. In the sequel we shall often take fibre products. We shall always assume that we are taking the reduced structure. When
we write that a diagram is ``cartesian modulo nilpotents", we mean that
the fibre product
is to be taken with the reduced structure. 
We shall always take care to define the cycle $[\Gamma]$.}

\end{rmk}

\subsection{Composition of correspondences}
\label{coc}
This formalism of correspondences via refined Gysin maps
extends easily to the case of composition of correspondences. The reference
is \ci{fulton}, Remark 16.1, Proposition 16.1.2,  $\S6$ and $\S 8$.

\bigskip
Let $X_i$, $i=1,2,3$ be nonsingular varieties. Let $p_{ij}:
X_1 \times X_2 \times X_3 \to X_i \times X_j$ be the obvious projections.
Let $\Phi \subseteq X_1 \times X_2$ and $\Gamma \subseteq X_2 \times X_3$
be two irreducible cycles such that all the following maps are proper:
$ |\Phi| \to X_2$, $\Gamma  \to X_3$, $p_{12}^{-1} (|\Phi|) \cap
p_{23}^{-1} (|\Gamma|) \to X_1 \times X_3$, $p_{13}(p_{12}^{-1} (|\Phi|) \cap
p_{23}^{-1} (|\Gamma|))  \to X_3$.
In addition, assume that 
$p_{13}(p_{12}^{-1} (|\Phi|) \cap
p_{23}^{-1} (|\Gamma|))$ has the expected dimension.

\medskip
We can therefore define a refined cycle $\Gamma \circ \Phi$ in 
$Z_{\ast}(p_{13}(p_{12}^{-1} (|\Phi|) \cap
p_{23}^{-1} (|\Gamma|) )$ and maps: $\Phi_*$, $\Gamma_*$
and $(\Gamma \circ \Phi )_*$. 
We have
$\Gamma_* \circ \Phi_* =  (\Gamma \circ \Phi )_*$.

\subsection{Quotient varieties}
\label{qv}
The formalism recalled so far
extends to the case of quotient varieties. See \ci{fulton}, Ex. 16.1.13.
In this case rational coefficients are necessary.
The basic set-up we need is as follows.

$Z$ is a nonsingular integral variety; $Z'$ is an integral variety; $p_Z: Z' \to Z$
is a surjective morphism; $G$ is a finite group acting on $Z$ and $Z'$, compatibly with
$p_Z$. There is  a commutative diagram:

\begin{diagram}
Z'           &   \rTo^{p'}           &  Z_1'  \\
\dTo^{ p_Z}  &                       &  \dTo^{p_{Z_1}}  \\
Z            & \rTo^{p}   &  Z_1\\
\end{diagram}

where $Z_1=Z/G$, 
$p$ is a quotient Galois  morphism of degree $|G|$, $Z_1'=Z'/G$ and $p'$
is a quotient Galois morphism of degree $|G|$. Note that it follows that
$p$ and $p'$ are separable.
Define $[Z']$ and $[Z'_1]$ to be the corresponding fundamental classes.

\bigskip

We have that
${p'}^*[Z_1']=[Z']$.
We can define a Gysin-type map as follows: 
$$
\Phi_{[Z_1']} : A_\ast(Z_1) \longrightarrow  A_\ast(Z'_1)\, , \qquad \quad
z_1  \longrightarrow    \,\frac{1}{|G|}\, p'_*(    p_Z^*(p^* z_1) \cap [Z']). 
$$
%If we fix $z_1'$, then we get a map
%$$
%A(Z_1) \to A({Z_1'}) \qquad z_1 \to z_1\cdot z_1'
%$$
%We have the factorization:
%$$
%A(Z_1) \stackrel{\frac{1}{|G|}p^*}\to A(Z) \stackrel{- \cdot_{Id_Z}
% {p'}^* z_1'}\to A(Z') \stackrel{{p'}_*}\to A(Z_1').
%$$
We need the following elementary fact which can be checked directly using the definition
of $p^*$ for quotient maps; see \ci{fulton}, Ex. 1.7.6. 
\begin{lm}
\label{symm}

Let $R$ be a variety endowed with an action of a finite group $G$ 
on it. Denote the quotient $p:R \to S$.  Let 
$\s: A_\ast(R) \to A_\ast(R)$, $ z \to \sum_{\gamma \in G}{\gamma_* z} $, be the so-called
symmetrization operator. 
Then
$
\s =  p^* p_*.
$
\end{lm}
%{\em Proof.} This can be checked directly at the level of cycles.
%We first assume that the characteristic of the base field is zero.

%Let $W \subseteq R$ be an integral cycle in $R$.  Let $F \subseteq
%G$ be the subgroup of $G$ fixing the cycle $W$ (i.e. $g(W)=W$). Let
%$I\subseteq F$ be the subgroup of elements which act as
%the identity of $W$.  Denote by $V$ the cycle
%$p(W) \subseteq W$. We have $p_*W = \deg{(W/V)}\,V$. Note that
%$\deg{(W/V)}= [F:I]$, the index of  $I$  in $F$.

%Let $W_i$ be the irreducible components of $f^{-1}(V)$. They are all isomorphic to %each other.
%In particular they are all isomorphic to W, which is one of them. The $I_i$
%all have the same cardinality. Same for the $F_i$. We thus have

%$ p^*p_*(W)= [F:I] p^*V= [F:I] \sum_{i}|I_i| W_i=  \sum_{i}{[F_i:I_i] |I_i| W_i}=
 %\sum_{i}{|F_i| W_i}=\s(W)$ ($W_i$ appears as many times as there are the elements
%of $F_i$).

%If the characteristic of the base field is positive, the morphism $p$ is still %separable, but
%a morphism $W \to V$ as above may fail to be separable. In this case the degree
%$\deg{(W/V)} = \deg{(W/V)}_s \cdot \deg{W/V}_i$ (separable and inseparable %degrees). Clearly,
%the general fiber of $W \to V$ has cardinality $\deg{(W/V)}_s$ $= [F:I]$.
%Since, by definition,  $p^*$ contains the inseparable degrees as denominators, the %same proof
%given above goes through.
%\blacksquare

\begin{lm}
\label{factoriz}

\medskip
We  have the following commutative diagram:

\begin{diagram}
A_{\ast}(Z')      &  \rTo^{p'_*}     &  A_{\ast}(Z_1')   \\
\uTo^{\Phi_{[Z']}}   &                 &  \uTo^{\Phi_{[Z_1']}}   \\
A_{\ast}(Z)       &  \rTo^{p_*}   & A_{\ast}( Z_1) \\
\end{diagram}

\end{lm}
{\em Proof.} 
For ease of notation we denote a 
class
$p_Z^*(a) \cap p'^*z'_1$, by $a\cdot p'^*z_1'$.

\n
Let $z \in A_\ast(Z)$ and $z' \in A_\ast(Z')$. Define
$\s: z \to \sum_{\gamma \in G}{\gamma_* z} $ and $\s': z' \to \sum_{\gamma \in G}{\gamma_* z'} $. 
By virtue of Lemma \ref{symm}, we have that
$
\s = p^* p_*$ 
and that 
$ \s'={p'}^* {p'}_*  .
$
It follows that 
$
{p'}^* {p'}_* (z \cdot {p'}^*z_1')= \s'(z \cdot {p'}^*z_1')= \sum_{\gamma \in G}{\gamma_*(
z \cdot {p'}^*z_1')}=$ (by the proper push-forward property) $=
\sum_{\gamma \in G}{(\gamma_*z \cdot \gamma_* {p'}^*z_1')}=$ (since  ${p'}^*z_1'$ is $G$-invariant)
$=\sum_{\gamma \in G}{(\gamma_*z \cdot  {p'}^*z_1')}$
$ =\s (z) \cdot  {p'}^* z_1=$ $(p^* p_* z) \cdot {p'}^*z_1'=$ (since 
$p^*p_* z $ is $G$-invariant) $= \frac{1}{|G|}
[\sum_{\gamma \in G}{\gamma_*(p^* p_* z ]\cdot {p'}^*z_1'   }=$
$\frac{1}{|G|} \s'(p^* p_* z)\cdot p'^* z_1'$ 
$
= \frac{1}{|G|}{p'}^* {p'}_* ( p^* p_* z \cdot {p'}^*z_1' ) $. 
Since ${p'}^*$ is injective, we have that 
$ {p'}_* (z \cdot {p'}^*z_1')=  \frac{1}{|G|}
{p'}_* ( p^* p_* z \cdot {p'}^*z_1' )$ which
is what we wanted to prove.
\blacksquare

\subsection{The morphisms $\Gamma_*$ and $\hat{\Gamma}_*$}
\label{mgg}
Let $X$ be an irreducible nonsingular quasi projective surface
defined over an algebraically closed field.

Recalling the notation introduced in $\S$\ref{cubo}, we  summarize the results of this section by stating that the following
diagram is commutative:

\begin{diagram}
A_{\ast}(\Gamma)= \oplus_{\nu} A_{\ast}(\Gamma^{\nu}) & \rTo & A_{\ast}(\hat{\Gamma})=\oplus_{\nu} A_{\ast}(\hat{\Gamma}^{\nu})  &   \rTo & A_{\ast}(X^{[n]})\\
   \uTo_{\Gamma_*=\oplus_{\nu} \Gamma_*^{\nu}}   &      &    \uTo_{\hat{\Gamma}_*=\oplus_{\nu} \hat{\Gamma}_*^{\nu}}            &        &                    \\
A_{\ast}({\cal X}) =  \oplus_{\nu} A_{\ast}(X^{\nu})      & \rTo & A_{\ast}(\hat{\cal X}) =   \oplus_{\nu} A_{\ast}(X^{(\nu)}).          &        &                     \\
\end{diagram}

By abuse of notation, the corresponding maps into $A_{\ast}(X^{[n]})$ will be denoted with the same symbol, e.g. $\Gamma_*: A_{\ast}({\cal X}) \to A(X^{[n]})$
and $\hat{\Gamma}_*: A_{\ast}(\hat{\cal X}) \to A(X^{[n]})$.

 \section{The Chow groups of $X^{[n]}$}
 Let  $X$ be   an irreducible nonsingular quasi projective surface
defined over an algebraically closed field. The reader should keep in mind $\S$\ref{cubo}.
\subsection{The injectivity of $\hat{\Gamma}_*$}
\label{injectivity}
We shall  freely use the formalism developed for correspondences over quotient varieties
in $\S$\ref{qv}.

\begin{lm}
\label{zeroifsupport}
Let $Z$ be  an irreducible algebraic scheme of dimension $n$, quotient of a nonsingular algebraic scheme via a finite group.
Let $V$ and $W$ be pure-dimensional cycles on $Z$  of dimensions $k$ and $l$.
Let $f:Z \to Y$ be a proper morphism of algebraic schemes.

\n
If $\dim{f(|V|\cap |W|)} < k+l-n$, then $f_*(V\cdot W)=0$.
\end{lm}
{\em Proof.} There is a well-defined refinement of the product
in $A_{k+l-n}(|V|\cap |W|)$ which must map to zero. See
\ci{fulton}, Ex. 8.3.12.
\blacksquare

\bigskip
The following lemma is essentially a reformulation of \ci{ellstrom}. Before stating it, we introduce some notation:
Let $\nu$ be a partition of $n$. For $\underline {x}\in X^{(\nu)}_{reg}$ we set $F_{\nu}=( \pi^{-1}( \underline{x} ) )_{red}$, which we identify 
via $P$ with $p_{\nu}^{-1}(\underline{x})$ (cf. \ref{identify}). We set $m_{\nu} : = (-1)^{n - l(\nu)} \prod_{j=1}^{l(\nu)}{\nu_j}$.
We denote by  $ [F_{\nu}]\cdot X^{[n]}_{\nu} \in A_0(F_{\nu})$ the refined intersection defined by the closed imbeddings of $F_{\nu}$ and  $X^{[n]}_{\nu}$
in $ X^{[n]}_{\geq \nu}$.
\begin{lm}
\label{ellings}
$\deg{([F_{\nu}]\cdot X^{[n]}_{\nu} )} = m_{\nu}$.
\end{lm}
{\em Proof.} The case $\nu=n^1$ is precisely the main result in \ci{ellstrom}.  
The general case  follows from  the K\"unneth formula.
\blacksquare

\bigskip
In the following two propositions we compute the compositions 
$^t\hat\Gamma^{\bnu} \circ \hat\Gamma^{\bmu}$. Here $^t-$ denotes, as
usual, the transposed correspondence; see \ci{fulton}, $\S16.1$.

\begin{pr}
\label{ab}
Let $\bmu \neq \bnu$ be two distinct partitions of $n$.
Then $\,^t\hat\Gamma^{\bnu}  \circ \hat\Gamma^{\bmu}  =0$ in 
$A_{l(\mu)+l(\nu)}(X^{(\bmu)}\times X^{(\bnu)}).$
\end{pr}
{\em Proof.} Consider $X^{(\bmu)}\times X^{[n]} \times X^{(\bnu)}$
together with the natural projections $p,$ $\pi$ and $q$ to
$X^{(\bmu)}\times X^{[n]}$, $X^{(\bmu)}\times  X^{(\bnu)}$ and
$ X^{[n]} \times X^{(\bnu)}$, respectively. By virtue of Lemma
\ref{zeroifsupport}, it is enough to show that
$\dim{\pi( p^{-1} \hat\Gamma^{\mu} \cap q^{-1}\,^t\hat\Gamma^{\nu}    )} < 
l({\bmu}) + l(\bnu)$. 

One checks directly
that  $\dim{\pi( p^{-1} \hat\Gamma^{\mu} \cap q^{-1}\,^t\hat\Gamma^{\nu}    )}=$
$\dim{ ( \overline X^{(n)}_{\bmu} \cap \overline X^{(n)}_{\bmu} ) }$. Since
$\bmu\neq \bnu$,
this last dimension
is strictly less than $2 \min{(l({\mu}), l({\nu}))}$.
\blacksquare

\begin{pr}
\label{aa}
$^t\hat\Gamma^{\nu}  \circ \hat\Gamma^{\nu}  =
m_{\nu}
\Delta_{X^{(\nu)}}$ in $A_{2l(\nu)}(X^{(\nu)}\times X^{(\nu)}).$
\end{pr}
{\em Proof.} 
The cycle $^t\hat\Gamma^{\nu}  \circ \hat\Gamma^{\nu}$ is supported on the diagonal of 
$X^{(\nu)}\times X^{(\nu)}$ which is irreducible of the expected dimension, therefore
$^t\hat\Gamma^{\nu}  \circ \hat\Gamma^{\nu}=c \, \Delta_{X^{(\nu)}}$, with $c \in \rat$. Let $\underline {x} \in X^{(\nu)}_{reg}$.
Since the  map $\hat\Gamma^{\nu}_{reg} \to X^{(\nu)}_{reg}$
is flat,  we have that  $\hat \Gamma^{\nu}_{*}([\underline {x}])=[F_{\nu}]$.
By virtue of Lemma  \ref{ellings} and  of the projection formula we have that
$$
c=\deg{ ( ^t\hat\Gamma^{\nu}_*  \circ \hat\Gamma^{\nu}_*([\underline{x}]) ) } = 
\deg \, ( ^t\hat\Gamma^{\nu}_* ([F_{\nu}]) ) = \deg{ ([F_{\nu}] \cdot X^{[n]}_{\nu} ) } =m_{\nu}.
$$
\blacksquare

\begin{cor}
\label{injj}
The natural map
$$
{\hat{\Gamma}}_*= \bigoplus_{\nu \in \prt}{ \hat{\Gamma}^{\nu}_* }:
A_{\ast} ( \hat{\cal X} ) = \bigoplus_{ \nu \in \prt}{A_{\ast}(X^{(\nu)}) } \longrightarrow
A_{\ast} ( X^{[n]} )
$$
is injective.
\end{cor}
{\em Proof.} Propositions \ref{aa} and \ref{ab} imply that 
$\, \!^t {\hat{\Gamma}}_* \circ \hat{\Gamma}_*$ is a non zero multiple
of the identity being the map induced by a non zero multiple of the diagonal.
\blacksquare.

\begin{rmk}
\label{injopen} {\rm An argument identical to the one given above shows that
 Corollary \ref{injj} holds if we replace all spaces by the ones obtained by base change
with respect to any open immersion $U \to X^{(n)}$.}
\end{rmk}

\subsection{A surjectivity statement}
\label{surjectivity}
The following
surjectivity statement  is essential to proving
Proposition \ref{mainsurj} which, in turn, constitutes the
surjectivity part of the main result of this paper,
Theorem \ref{maintm}.
The proof  will be given at the end of this section, after a series of preparatory
results.
\begin{pr}
\label{mainlemma}
Let $\mu \in \prt$: the map 
$P_* (p_{\geq \mu}^!(- )\cap [\G{\geq \mu}]): \Ch{\ast}{ {\cal X}_{\mu} } \to \Ch{\ast}{\Hil{n}{\mu} }$ is surjective.
\end{pr}

\bigskip
In order to prove Proposition \ref{mainlemma} we  need to make explicit 
the combinatorics involved.
Let $l$ be an integer, and $\RR_{l}$ be the set of decompositions 
$\r=\{I_1,\cdots,I_r\}$ of the set $\{1, \cdots, l\}$ into disjoint subsets,  
i.e. $\{1, \cdots, l\}= I_1 \coprod \cdots \coprod I_r$. 
For $\r\in \RR_{l}$ define ${\cal D}_\r = \{ (x_1,\cdots,x_l) \in \sn{l} \,:\, x_i = x_j \Leftrightarrow \, i,j\in I_k\, \hbox { for some } k \}$. 
For $\r=\{1\},\cdots \{l\}$ we have ${\cal D}_\r=\snn{l}{reg}= X^l\setminus Diagonals$. For $\r=\{1, \cdots ,l\}$, we have  that ${\cal D}_\r$ is the small
diagonal.

Let $ \PP_{\preceq \nu}(n)= \{ \, \mu \in \prt \, \hbox{ such that } \mu \preceq \nu \}$. 
Define 
$Q_{\nu} : \RR_{l(\nu)} \to  \PP_{ \preceq \nu}(n)$ 
by setting
$Q_{\nu}(\{I_1,\cdots,I_r\})=
(\sum_{i \in I_1} \nu_i, \cdots ,\sum_{i \in I_{l}} \nu_i)$.
Note that this partition is not ordered, i.e.
$\sum_{i \in I_1} \nu_i \not\geq \sum_{i \in I_2}\nu_i
                               \cdots \not\geq\sum_{i \in I_{l}} \nu_i$, but
see Remark \ref{order}.

\medskip
Let  $\nu : \sn {\nu} \to \overline \Sym{n}{\nu}$ be the already mentioned map defined by
$\nu(x_1,\cdots ,x_{l(\nu)})=\nu_1 x_1 + \cdots +\nu_{l(\nu)}x_{l(\nu)}$.
Note that $\nu^{-1}(\Sym{n}{\nu})=\snn{\nu}{reg}$. Note that $l(Q_{\nu} ( \{ I_1, \ldots, I_r  \}))=r$.
The proof of the following lemma is immediate:
\begin{lm}
\label{strata}
$$\sn {\nu} \times_{\Sy{n}} \Sym{n}{\mu}\neq \emptyset \Leftrightarrow \nu \succeq \mu \hbox{ and }
(\sn{\nu} \times_{\Sy{n}} \Sym{n}{\mu})_{red}=\coprod_{\r \in Q^{-1}_{\nu}(\mu)}{\cal D}_\r . $$
\end{lm}

Let  $\Gn{\nu}{}$ be the correspondences introduced in  $\S$\ref{definitions}.

\begin{rmk}
\label{azione}
{\rm The group $\Si_{\nu}$ acts on $\RR_{l(\nu)}$ and $Q_{\nu}$ is 
$\Si_{\nu}$-invariant. In particular, $\Si_{\nu}$ acts on 
$Q_{\nu}^{-1}(\mu)$.}
\end{rmk}

\begin{rmk}
\label{order}
{\rm For $\nu \in \prt$ we can fix a total order on the set $\RR_{l(\nu)}$ 
%subsets of $\{ 1,\cdots,l(\nu)\}$
so that we always have 
$\sum_{i \in I_1} \nu_i \geq \sum_{i \in I_2}\nu_i
                               \cdots \geq\sum_{i \in I_{r}} \nu_i$.
                               
Let $\rho \in \RR_{l(\nu)}$ and $\mu:= Q_{\nu}(\rho)$. 
We can  identify ${\cal D}_{\r}$ with $\snn{\mu}{reg} $ 
via the map $\delta_{\r}: \snn{\mu}{reg} \to \sn{\nu}$
sending the point 
$(y_1,\cdots,y_{r})\in \snn{\mu}{reg}$ to $(x_1,\cdots,x_{l(\nu)})$ defined by $x_j=y_i$ 
precisely if $j \in I_i$. 
In this case $\nu \circ \delta_{\r}=\mu_{|\snn{\mu}{reg}}:\snn{\mu}{reg} \to \sn{(n)}$,
hence 
$(\Gn{\nu}{}\times_{\sn{\nu}}{\cal D}_{\r})_{red}
=
\Gn{\mu}{|\snn{\mu}{reg}}=:\Gn{\mu}{reg}$. 
These identifications will 
always be tacitly made in the sequel.}
\end{rmk}

\begin{rmk}
\label{fibre}
{\rm $p_{\mu}:\Gn{\mu}{reg}\to \snn{\mu}{reg}$ is a Zariski locally trivial fibration. 
The fiber is isomorphic to the product of punctual Hilbert schemes
${\cal H}_{\mu_1} \times \cdots \times {\cal H}_{\mu_l}$ and it admits a cellular decomposition. See \ci{go}, Lemma 2.1.4. and \ci{e-s1}.
}
\end{rmk}

%In view of remark \ref{order} every couple 
%$(\nu,\r),\hbox { with } \nu \in \prt,\hbox{ and }\r \in \RR_{l(\nu)}$ 
%such that $Q(\r)=\mu$ gives a diagram

%\begin{diagram}
%\Gn{\mu}{reg}         & \rTo         & \Gn{\nu}{}             \\
%\dTo^{p_{\mu}}       &              & \dTo^{p_{\nu}}      \\
%\snn{\mu}{reg}   & \rTo^{\delta_\r} & \sn{\nu}         \\
%\end{diagram}

Let  $(\nu,\r)$ be such that $Q_{\nu}(\r)=\mu$.  The pair $(\nu,\r)$ defines  a set of partitions 
$\nu^i \in \PP(\mu_i)$ as follows: if $\r=\{I_1,\cdots,I_{l(\mu)}\}$, then    $\mu_i= \sum_{j \in I_i}\nu_j$ and therefore 
$\nu^i=\{\nu_j\}_{j \in I_i}$ is a partition of $\mu_i$.

\begin{defi}
Define the open sets 
%$U_{\r} \subseteq \sn{\nu}$ by 
$$
U_\r=\{ (x_1,\cdots,x_{l(\nu)}) \in X^{\nu}, \hbox { such that the subsets } 
                                       \{ x_j\}_{j\in I_i} \hbox{ are pairwise disjoint} \}.$$                        
\end{defi}  

In other words we allow only points belonging to the same $I_j$ to collapse:
in particular $U_{\r} \supseteq {\cal D}_{\r}$.

\begin{lm}
\label{split}
 Let  $(\nu,\r)$ be such that $Q_{\nu}(\r)=\mu$, and $\nu^i \in \PP(\mu_i)$ 
be the corresponding set of partitions. 
There is a canonical isomorphism 
$$
\Gn{\nu}{|U_{\r}}=\prod \Gn{\nu^i}{|U_{\r}}.
$$

\end{lm}

{\em Proof.}
Let $( x_1, \cdots,x_{l(\mu)}, {\cal I} ) \in \Gn{\nu}{|U_{\r}}$. From the definition
of $U_{\r}$ it follows that ${\cal I}$ is the product of ideals 
${\cal I}_1,{\cal I}_2,\cdots, {\cal I}_r$ of lengths $\mu_1, \mu_2,\cdots,\mu_r$ supported on
 $\{ x_j\}_{j\in I_1},  \{ x_j\}_{j\in I_2}, \ldots ,$
$  \{ x_j\}_{j\in I_r}$, whence the obviously bijective map
from $\Gn{\nu}{|U_{\r}}$ to $\prod \Gn{\nu_i}{|U_{\r}}$.
\blacksquare

\bigskip
We have the following  diagram 
\begin{diagram}
\Gn{\mu}{reg}                & \rTo              & \Gn{\nu}{|U_{\r}}    \\
\dTo^{p_{\mu}}               &                   & \dTo^{p_{\nu}}       \\
{\cal D}_{\r}= \snn{\mu}{reg}  & \rTo^{\delta_\r}  & U_{\r}    .           \\
\end{diagram}
which is cartesian modulo nilpotents (see Remark \ref{cmn}) and a map 
$\gamma_{\nu,\r}: \Ch{\ast}{\snn{\mu}{reg}}\to \Ch{\ast}{\Gn{\mu}{reg}}$
   defined by 
$\gamma_{\nu,\r}(\beta)=p_{\nu}^!(\beta)\cap [\Gn{\nu}{}]$ (cf. $\S$\ref{intersection}).

\begin{lm}
\label{pullback}
Let $\alpha_{\nu,\r}:=\gamma_{\nu,\r}([\snn{\mu}{reg}]) \in \Ch{\ast}{\Gn{\mu}{reg}}$
and $\beta \in \Ch{\ast}{\snn{\mu}{reg}}$. Then 
$$\gamma_{\nu,\r}(\beta) = p_{\mu}^*(\beta)\cap \alpha_{\nu,\r} . $$
\end{lm}
{\em Proof.}
It follows from  the associativity
property of refined products, \ci{fulton}, Proposition 8.1.1.a).
\blacksquare

\bigskip
What has been done so far can be summarized by the following:
\begin{lm}
\label{summary}
Let $\mu$ be a partition, denote by 
$\tilde q: \Gamma_{reg}^{\mu} \to \Hil{n}{\mu}$  
the quotient map by the action of $\Si_{\mu}$. Let $W \subseteq A_{\ast} (\Gamma_{reg}^{\mu})$ be the 
$\Ch{\ast}{\snn{\mu}{reg}}$-submodule  generated by  the classes
$\{ \alpha_{\nu , \r } \}_{\r \in Q_{\nu}^{-1}(\mu)}$.
The surjectivity of the map 
$P_* (p_{\geq \mu}^!( - )\cap [\G{\geq \mu}]): 
\Ch{\ast}{ {\cal X}_{\mu} } \to \Ch{\ast}{\Hil{n}{\mu} }$ is 
equivalent to the surjectivity of the restriction of
$\tilde q_*: A_{\ast}(\Gamma_{reg}^{\mu}) \to 
A_{\ast}(\Hil{n}{\mu})$ to $W$.  
\end{lm}
{\em Proof.} By virtue of  lemma \ref{strata} 
$$
{\cal X}_{\mu}=\coprod_{\r \in Q_{\nu}^{-1}(\mu)} {\cal D}_{\r}
= \coprod_{\r \in Q_{\nu}^{-1}(\mu)} X^{\mu}_{reg},
$$
this last identification is made using   the map described in \ref{order}.                                               
The following diagram, where the equality 
$\sum \gamma_{\nu,\r}=\sum p_{\mu}^*(\,)\cap \alpha_{\nu, \r}$
is a consequence of \ref{pullback}, commutes:
\begin{diagram}
\Ch{\ast}{ {\cal X}_{\mu} } & \rTo^{ P_* (p_{\geq \mu}^!(\,)\cap [\G{\geq \mu}])} & \Ch{\ast}{\Hil{n}{\mu} } \\
\dTo                           &                                                     & \uTo^{\tilde q_*}           \\
\oplus_{\r \in Q_{\nu}^{-1}(\mu)} A_{\ast}(X^{\mu}_{reg})& \rTo^{\sum \gamma_{\nu,\r}=\sum p_{\mu}^*(\,)\cap \alpha_{\nu, \r}}& A_{\ast}(\Gamma^{\mu}_{reg}) ,\\
\end{diagram}
whence the statement.
\blacksquare

\bigskip                                                                                                       
We  now study a special case which is crucial to the proof of Proposition
\ref{mainlemma}.
We consider $\mu=n^1$.  
For every partition $\nu$ there is only one  
$\delta_{\r}= \delta : X \to \sn{l(\nu)}$ and $p:= p_{n^1}: \Gn{n^1}{} \to X $ is a Zariski locally trivial fibration, 
whose fibre is isomorphic to the length $n$ punctual Hilbert scheme ${\cal H}_n$.
For each $\nu \in \prt$ we have $\alpha_{\nu}=
                   p_{\nu}^!([X])\cap [\Gn{n^1}{}]$ 
and the map
$P: \Ch{\ast}{X}^{\oplus p(n)} \to \Ch{\ast}{\Gn{n^1}{}}$ sending the collection $\{\beta_{\nu}\}$ 
to the class $\sum_{\nu} p^{*}(\beta_{\nu})\cap \alpha_{\nu}$.
\begin{lm}
\label{casobase}
The map $P$ is surjective
\end{lm}
{\em Proof.}
By  virtue of Lemma \ref{fibration}, it is enough to prove that the restrictions $\hat \alpha_{\nu}$ of the classes
$\alpha_{\nu}$ to any fibre generate $\Ch{\ast}{{\cal H}_n}$.

We may assume, without loss of generality that $X$ is projective, for the restrictions $\hat{\alpha}_{\nu}$ do not change.
Let $i_x: \{point\} \to X$ be the imbedding of a point $x$. Note that 
$\hat \alpha_{\nu}=i_x^!([\Gn{\nu}{}])$. Let $g: {{\cal H}_n}  \to
X^{[n]}$ be the closed embedding and define
$\hat{\beta}_{\nu} : =g_*\hat \alpha_{\nu}$ 
$ \in A_{\ast}(X^{[n]})$. 

Consider the pairing $A_{\ast}(X^{[n]})  \times A_{\ast}(X^{[n]}) \to \rat$
given by taking $(a,b) \to \deg{a \cdot b}= \int_{ X^{[n]} }{ a\cdot b}$, where the last product is the product in the Chow ring. Note that this pairing is almost never perfect. It descends to
algebraic equivalence.

\medskip
CLAIM. 
$\deg{ ( \hat \beta_{\nu} \cdot [\overline \Hil{n}{\mu}] ) }=0$ if and only if 
$\mu \neq \nu$.

Proof of the CLAIM. Consider $\Gamma^{\nu}$  and $X^{\nu}\times \overline X^{[n]}_{\nu}  $ as a pair of  cycles
in $X^{\nu} \times X^{[n]}$. They  define two families of cycles on $X^{[n]}$:
$\Gamma^{\nu}_{\tau}$ and 
$\{ {X^{\nu}\times \overline X^{[n]}_{\nu} } \}_{\tau} $, $\tau \in  X^{\nu}$. See \ci{fulton}, $\S10$. 
In particular,  we have
$ \Gamma_y^{\nu}= i_y^!([\Gamma^{\nu}])$.
If $y\in \snn{\nu}{reg}$, then 
$\Gamma^{\nu}_y=[{(p_{\nu}^{-1}(y)})_{red}]$. 
Recalling  that  we have   the canonical $\delta_{\rho}: X \to X^{\nu}$,
we see that
$\hat\beta_{\nu}= \Gamma^{\nu}_x$. 
In other words: $ \Gamma_y^{\nu}$ and $\hat\beta_{\nu}$ belong to the same family of cycles.
 By construction, $\hat\beta_{\nu}$
can be represented, modulo algebraic equivalence, 
by a cycle supported at ${( p_{\nu}^{-1}(y) ) }_{red}$.
 
Note that  $\hat \beta_{\nu} \cdot  [\overline \Hil{n}{\mu}] \in A_{l(\nu) - l(\mu)}(X^{[n]})$.

If $l(\nu) < l(\mu)$, or if $l(\nu) > l (\mu)$, then the degree is zero for trivial reasons.

If $l(\nu) = l(\mu)$, but $\mu \neq \nu$, then the degree is still zero. 
In fact, using refined intersection products in the context of algebraic equivalence, we
can represent $\hat\beta_{\nu}$ as a cycle supported on a fiber $\pi^{-1}(y)$. This fiber  
 does not meet
$\overline \Hil{n}{\mu}$.

Finally, if $\nu =\mu$, then, by virtue of \cite{fulton},
Corollary 10.1 and Proposition 10.2, we have that 
$$
\deg{ ( \hat{\beta}_{\nu} \cdot [\overline\Hil{n}{\mu}] ) }=
\deg{ ( \Gamma_y \cdot  [\overline\Hil{n}{\mu}]  ) } = 
  \deg ( [p_{\nu}^{-1}(y)] \cdot [\overline X^{[n]}_{\nu}])=m_{\nu}.
$$

The CLAIM follows easily.

By virtue of Corollary \ref{injj} (see also
\ci{decatmig}), 
the classes $[\overline \Hil{n}{\mu}]$ are independent.
Jointly with  the CLAIM just proved,
this shows that the classes $\hat{\beta}_{\nu}$  are independent.
It follows that the classes   
  $\hat{\alpha}_{\nu}$ are independent.
Their number is $p(n)$,
which is equal to the dimension of $\Ch{\ast}{{\cal H}_n}$ (cf. \ci{decatmig}, for example).
 This proves the 
lemma.
\blacksquare

\medskip

Let us go back to the situation dealt with in Lemma \ref{split}: the 
diagram after that lemma and \ref{casobase} 
imply immediately:

\begin{lm}
\label{classes}
The restriction of the cycle $\alpha_{\nu,\r}:=p_{\nu}^!([{\cal D}_{\r}])\cap [\Gn{\nu}{}] \in \Ch{\ast}{\Gn{\mu}{}}$ 
to a fibre of $p_{\mu}$ is the cycle
$\hat \alpha_{\nu^1} \times \cdots \times \hat \alpha_{\nu^r} 
                                   \in \Ch{\ast}{{\cal H}_{\mu_1} \times \cdots \times {\cal H}_{\mu_l}}$.
\end{lm}

Let $\underline \eta=(\eta^1, \cdots,\eta^r)$ be a multipartition of $\mu$. 
By this we mean that  $\eta^i \in \PP(\mu_i)$.
Set $l_i=l(\eta^i)$ and define the open sets 
$U_{\underline \eta} \subseteq \prod \sn{\eta^i}$ by 
$$
U_{\underline \eta}= \{ (x_1^{(1)},\cdots ,x_{l_1}^{(1)}, \cdots , x_1^{(r)},\cdots,x_{l_r}^{(r)} ) 
                                        \;   \hbox{ such that } \;  x_k^{(i)} \neq x_l^{(j)} \hbox{ if } i \neq j \} .
$$
In other words we allow collisions only among points of the same group.
$\snn{\mu}{reg}$ can be identified with the closed subset of points in $U_{\underline \eta}$ such that
 $x_k^{(i)} = x_l^{(j)} \hbox{ if and only if } i = j $.
 
According to \ci{fulton}, 8.1.4, the diagram

\begin{diagram}
\prod \Gn{\mu_i}{|\snn{\mu}{reg}}                       & \rTo         & \prod \Gn{\eta^i}{}                 \\
\dTo^{\times p_{\mu_i}}                                 &              & \dTo^{\times p_{\eta^i}}              \\
\snn{\mu}{reg}                                          & \rTo         & \prod \sn{\eta^i}                      \\
\end{diagram}

produces classes $\alpha_{\underline \eta }\in \Ch{\ast}{\prod \Gn{\mu_i}{|\snn{\mu}{reg}}}$,
$\alpha_{\underline \eta }=(\times p_{\eta^i})^!([\snn{\mu}{reg} ])\cap [ \prod \Gn{\eta^i}{} ]$ 
whose restrictions to a fibre is the cycle
$\hat \alpha_{\eta^1} \times \cdots \times \hat \alpha_{\eta^r} 
\in \Ch{\ast}{{\cal H}_{\mu_1} \times \cdots \times {\cal H}_{\mu_l}}.$

\begin{lm}
\label{span}
The set  $\{ \alpha_{\underline \eta } \}$ generates $ \Ch{\ast}{\Gn{\mu}{reg}}$ as a  $\Ch{\ast}{\snn{\mu}{reg}}$-module.
\end{lm}
{\em Proof.} Since ${\cal H}_{\mu_1} \times \cdots \times {\cal H}_{\mu_r}$ has a cellular decomposition,
$\Ch{\ast}{{\cal H}_{\mu_1} \times \cdots \times {\cal H}_{\mu_r}}=\Ch{\ast}{{\cal H}_{\mu_1}} \otimes \cdots 
            \otimes \Ch{\ast}{{\cal H}_{\mu_r}}$ has dimension $\prod p(\mu_i)$.
The set in question is therefore a basis of this space, and the statement is  a direct consequence of Lemma \ref{fibration}, Remark \ref{fibre} and Lemma \ref{classes}.
\blacksquare

\begin{rmk}
{\rm We have already observed that a couple $(\nu,\r)$ with 
$Q_{\nu}(\r)=\mu$ 
determines a multipartition $\underline \nu_{\r}$ of $\mu$;
Lemma \ref{split} implies  that   
$\alpha_{\underline \nu_{\r}}=\alpha_{\nu,\r}$.}
\end{rmk}

\begin{rmk}
{\rm The set $\{\alpha_{\nu,\r}\}$ is only a subset of 
$\{ \alpha_{\underline \eta}\}$ 
as shown by the following example}
\end{rmk}

\begin{ex}
{\rm Let $n=4$, and $\mu=2^2$. There are $4$ multipartitions: $(2,2)$, $(2,1^2)$, $(1^2,2)$,$(1^2,1^2)$
and $4$ corresponding cycles which restrict to a basis for $\Ch{\ast}{{\cal H}_2 \times {\cal H}_2}$.
The couples $(\nu, \r)$  such that $Q_{\nu}(\r)=2^2$ are:
if $\nu=1^4$, then  there are $3$ different $\r$'s which are conjugate by $\Si_{\nu}$ and give the same (0-dimensional) cycle
corresponding to the multipartition $(1^2,1^2)$;
if $\nu=2^2$, then  there is only one $\r$ which gives the 2-dimensional cycle corresponding to $(2,2)$;
if $\nu=2 \cdot 1^2$, then  the only $\r$ is $\{1\},\{2,3\}$ which induces the multipartition $(2,1^2)$.
This is consistent with the fact that $\sn{ 2\cdot 1^2} \times_{\Sy{4}} \Sym{4}{2^2}$ contains only 
one component ${\cal D}_{\r}$ given by points of  type $x_2=x_3$. Note also that the map 
${\cal D}_{\r} \to \Sym{4}{2^2}$ is $2:1$, the quotient map by $\Si_{2^2}$. }
\end{ex}

This can be easily explained: note first that the group
$\Si_{\mu}$ acts on the set of multipartitions of $\mu$.
Given a multipartition 
$\underline \eta=(\eta^1, \cdots,\eta^r)$
of $\mu$,
let $\eta^i=\eta_1^i \geq \cdots \geq \eta_{l_i}^i$.
The sequence of the $\eta_j^i$ is a partition $\nu \in \prt$, let  $l= \sum l_i$ be its length. 
A permutation $\s \in \Si_l $, reordering the sequence in a non-increasing one, is identified up to left multiplication with $\Si_{\nu}$.
Define the subsets 
%\bigskip
%???????????????????????
%DIMMI SE QUESTO VA MEGLIO: PENSA A $\s$ come a una applicazione biiettiva da $\{1, \cdots,l\}$ in se'.
$I_i:=\{ \s((\sum_{k=0}^{i-1}l_k)+1),\cdots,\s(\sum_{k=0}^i l_k) \}.$
%{\bf so that after reordering  as in \ref{order} ???????? AGGIUSTI TU QUA?} 
We get 
$\{I_i\} \in \RR_{l} $, and $\sum_{k \in I_i}\nu_k=\mu_i$.
This associates with $\underline \eta $ a couple $(\nu,\r)$ 
such that $ \r \in Q_{\nu}^{-1}(\mu)$.
A different choice of the permutation $\s$ gives the same partition $\nu$
whereas $\r$ is changed by the action of $\Si_{\nu}$ (cf. Remark \ref{azione}).
Starting from the couple $(\nu,\r)$, instead, gives, as we have already observed, 
a multipartition $\hat \eta $ of $\mu$.
Note  that this multipartition depends on the way the $I_i$'s were ordered. 
We thus have:
\begin{lm}
$\hat \eta$ and $\underline \eta$ are in the same $\Si_{\mu}$-orbit.
\end{lm}

\medskip

Finally, let  $q:{\cal D}_{\r}=X^{\mu}_{reg} \to \Sym{n}{\mu}$ and
$\tilde q: \Gamma_{reg}^{\mu} \to \Hil{n}{\mu}$  be the quotient maps by $\Si_{\mu}$.
It follows from \cite{fulton}, Ex. 1.7.6,  that $q^*$ (resp. $\tilde q^*$) identify 
$\Ch{\ast}{\Sym{n}{\mu}}$ (resp. $\Ch{\ast}{\Hil{n}{\mu}}$) with 
$\Ch{\ast}{{\cal D}_{\r}}^{\Si_{\mu}}$, (resp. $\Ch{\ast}{\Gn{\mu}{reg}}^{\Si_{\mu}}$).
By virtue of lemma \ref{symm},  if 
$\alpha \in \Ch{\ast}{{\cal D}_{\r}}$, then $q^*q_* \alpha = \sum_{\s \in \Si_{\mu}} \s \cdot \alpha$ . Similarly,
if    $\alpha \in \Ch{\ast}{\Gn{\mu}{reg}}$, then 
$\tilde q^* \tilde q_* \alpha = \sum_{\s \in \Si_{\mu}} \s \cdot \alpha$.

\begin{lm}
\label{simm}
The images under the map $\tilde q^* \tilde q_*$ of the 
$\Ch{\ast}{\snn{\mu}{reg}}$-submodules  generated by
$\{ \alpha_{\underline \eta } \}$ and $\{ \alpha_{\nu , \r } \}$ coincide.
\end{lm}
{\em Proof.} We fix our attention on a single $\Si_{\mu}$-orbit. Set
$ \alpha_{\nu , \r }=\alpha_0$. Let $H$ denote the  stabilizer of $\alpha_0$ 
in $\Si_{\mu}$, and choose a set $g_1=e,g_2, \ldots ,g_r$ of representatives
for $\Si_{\mu}/H$, so that the orbit of $\alpha_0$ is
$\alpha_0,g_2 \alpha_0,\ldots,g_r \alpha_0$.
Given a cycle $\alpha=\sum_{i=1}^r p^*(\beta_i)g_i\alpha_0$, let
$\beta=(\sum_{i=1}^r g_i^{-1} p^*(\beta_i))\alpha_0$.
Then $\tilde q^* \tilde q_* \alpha=\tilde q^* \tilde q_* \beta$.
In fact, 
$
\tilde q^* \tilde q_* \alpha =\sum_{g,j}gp^*(\beta_j)gg_j\alpha_0
=
\sum_ig_i(\sum_{\stackrel{h\in H} {j=1}}^r hg_j^{-1}p^*(\beta_j))g_i\alpha_0 ,
$ 
while
$
\tilde q^* \tilde q_* \beta=\sum_g g(\sum_j g_j^{-1}p^*(\beta_j))g\alpha_0
=
\sum_i(\sum_{\stackrel{h \in H} {j=1}}^r g_ihg_j^{-1}p^*(\beta_j))g_i\alpha_0 .
$
\blacksquare   
 
\bigskip
We are now in the position to prove Proposition \ref{mainlemma}.  
Lemma \ref{span} and
Lemma \ref{simm} imply that
the restriction of the push-forward map $\tilde q_*$ to the 
$\Ch{\ast}{\snn{\mu}{reg}}$-submodule of 
$A_{\ast}(\Gamma^{\mu}_{reg})$ generated by the classes
$\{ \alpha_{\nu , \r } \}$ is surjective.
This implies Proposition \ref{mainlemma} by  virtue of Lemma
\ref{summary}.

\subsection{The  surjectivity of $\hat{\Gamma}_*$}

\label{cghs}
Proposition \ref{mainlemma} implies easily the following
surjectivity result:

\begin{pr}
\label{mainsurj}
Let $[\Gamma]$ be the fundamental cycle of $\Gamma$. 
The map 
$$
P_* (p^*( - )\cap [\Gamma])= \Gamma_*: \Ch{\ast}{ {\cal X}} =
\bigoplus_{\nu \in \prt}{ A_{\ast}( X^{\nu} ) }\to \Ch{\ast}{\Hi{n}}
$$
is surjective.
\end{pr}
{\em Proof.}
By virtue of Lemma \ref{exact} we have, for every $\mu \in \prt$, a commutative diagram with exact columns:
\begin{diagram}
\Ch{\ast}{{\cal X}_{\mu}}      & \rTo^{p_{\geq \mu}^!( - )\cap [\G{\geq \mu}]} & \Ch{\ast}{\G{\mu}} & \rTo^{P_*} & \Ch{\ast}{\Hil{n}{\mu}}                            \\
\dTo^{i_*}                      &               & \dTo                                       &                                      & \dTo                                        \\
\Ch{\ast}{{\cal X}_{\geq \mu}} & \rTo^{p_{\geq \mu}^*( -)\cap [\G{\geq \mu}]}           & \Ch{\ast}{\G{\geq \mu}}              & \rTo^{P_*} & \Ch{\ast}{\Hil{n}{\geq \mu}} \\
\dTo^{j^*}                      &                                    & \dTo                                          &              & \dTo  \\
\Ch{\ast}{{\cal X}_{> \mu}}    & \rTo^{p_{ > \mu}^*( - )\cap [\G{> \mu}]} & \Ch{\ast}{\G{> \mu}}      & \rTo^{P_*} & \Ch{\ast}{ \Hil{n}{ > \mu}} \\ 
\dTo         &             &   \dTo               &              & \dTo  \\
0             &             &    0                 &              &  0 .     \\              
\end{diagram}
We proceed by decreasing induction on $\mu$, proving that  
$P_* (p^*( - )\cap [\Gamma_{\geq \mu}]) : \Ch{\ast}{ {\cal X}_{\geq \mu}} \to \Ch{\ast}{\Hil{n}{\geq \mu}}$ 
is surjective. The statement in the case $\mu=1^n$ 
is clearly true since $\Hil{n}{1^n}=\Sym{n}{1^n}={\cal X}_{1^n}/ \Si_{n}$.
Suppose now that the surjectivity of $P_* (p^*( - )\cap [\Gamma_{>\mu}]): \Ch{\ast}{ {\cal X}_{>\mu}} \to \Ch{\ast}{\Hil{n}{> \mu}}$ has been established. By virtue of Proposition  \ref{mainlemma},  the map
in the first line of the diagram is surjective, whence the surjectivity of the second line and the statement.
\blacksquare
 
\begin{cor}
\label{surjquot} 
The natural map 
$$
\hat{\Gamma}_*: A_{\ast}(\hat{\cal X})=
\bigoplus_{\nu \in \prt}{ A_{\ast} ( X^{ (\nu) } ) }
 \longrightarrow  A_{\ast}(X^{[n]})
$$
is surjective.
\end{cor} 
{\em Proof.} Immediate from Proposition \ref{mainsurj} and $\S$\ref{mgg}.

\begin{rmk}
\label{surjopen}
{\rm Corollary \ref{surjquot} holds if we replace all spaces by the ones obtained by base change
with respect to any open immersion $U \to X^{(n)}$. In fact, it is sufficient to
use the corollary together with the standard exact sequence
\ci{fulton}, Proposition 1.8.}
\end{rmk}

\subsection{The isomorphism of Chow groups}
\label{icg}
The main result of this paper now follows easily. It does {\em not} 
hold with $\zed$ coefficients.
\begin{tm}
\label{maintm}
Let $X$ be an irreducible nonsingular algebraic surface defined over an algebraically closed
field.
The natural map
$$
\hat{\Gamma}_*= \bigoplus_{\nu \in \prt}{\hat{\Gamma}^{\nu}_*}: \; 
A_{\ast} (\hat{\cal X}) = \bigoplus_{\nu \in \prt}{A_{\ast}(X^{(\nu)})} \longrightarrow
 A_{\ast} (X^{[n]})
$$
is an isomorphism.
\end{tm}
{\em Proof.} Injectivity and surjectivity are proved in Corollary \ref{injj}
and \ref{surjquot}, respectively.
\blacksquare

\bigskip

The following is an immediate consequence
of Theorem \ref{maintm}. Let $Y$ be a scheme.
We denote by $K_o(Y)$ the Grothendieck
group of coherent sheaves on $Y$ and  we define
$K_o(Y)_{\rat}:= K_o(Y) \otimes_{\zed} {\rat}$.
\begin{cor}
\label{groth}
There is a  natural isomorphism of Grothendieck groups
with $\rat$-coeffi\-cien\-ts
$$
\hat{\Gamma}_*= \bigoplus_{\nu \in \prt}{\hat{\Gamma}^{\nu}_*}:
\; K_o (\hat{\cal X})_{\rat} = \bigoplus_{\nu \in \prt}{K_o(X^{(\nu )})_{\rat}} \longrightarrow
 K_o (X^{[n]})_{\rat}.
$$
\end{cor}
{\em Proof.} It follows immediately from Theorem 
\ref{maintm} and  \ci{fulton}, Corollary 18.3.2.
\blacksquare

\begin{rmk}
\label{equi}
{\rm 
Corollary \ref{groth} and the localization theorem for  equivariant
$K$-theory \ci{vi}, imply that one has an isomorphism
$K_o^{\Sigma_n}(X^n)_{\rat} \simeq K_o(X^{[n]})_{\rat}$. In a paper in preparation, we construct this isomorphism directly, i.e. without the localization theorem.
We expect similar statements to hold for the higher $K$-theory and for
the equivariant
derived category. 
}
\end{rmk}

\begin{rmk}
\label{obc}
{\rm Theorem \ref{maintm} and Corollary \ref{groth} hold
if we replace $\cal X$ and $X^{[n]}$ by the corresponding open sets
obtained by base change with respect to any open immersion $U \to X^{(n)}$.
See Remark \ref{injopen} and \ref{surjopen}.}
\end{rmk}

\begin{rmk}
{\rm Note that, by virtue of \ci{fulton}, Ex. 1.7.6,  $A_{\ast}(X^{(\nu)})= A_{\ast} (X^{\nu})^{\Sigma_{\nu}}$.
However, unless $X$ itself has a cellular decomposition, the natural map
$A_{\ast} (X)^{\otimes l(\nu)} \to A_{\ast} (X^{(\nu)})$ is not necessarily surjective.
 This contrasts singular cohomology, where one has
K\"unneth Formula.} 
\end{rmk}

\begin{rmk}
\label{anyf} 
{\rm With minor modifications, which we leave to the interested reader,
Theorem \ref{maintm}  holds: 1)  for a geometrically
irreducible smooth surface $X$ defined over any field, and  2)  in the analytic context
where $X$ is a smooth complex analytic surface (see \ci{fulton}, Ex. 19.2.5). }
\end{rmk}

\section{The motive of $X^{[n]}$ }
\label{motive}
Let $X$ be an irreducible nonsingular quasi projective surface 
defined over an algebraically closed field.

\subsection{The correspondences $\Delta_{\nu}$}
\label{cdn}
We freely use refined products over quotient varieties and the formalism of correspondences between
quotient varieties together with their standard properties. See
$\S$\ref{correspondences}.
Let $\nu \in \prt$ and consider $(X^{[n]}_{\nu}\times_{X^{(n)}}X^{[n]}_{\nu})_{red}$. Being 
the quotient of a Zariski locally trivial fibration with irreducible fiber
by the action of a finite group, this space  is a $2n$-dimensional irreducible locally closed subset of
$X^{[n]} \times X^{[n]}$. Let $D^{\nu}\subseteq X^{[n]} \times X^{[n]}$ be its closure:
$$
D^{\nu} := \overline{\{
(a,b) \in X^{[n]} \times X^{[n]} \, | \; \pi(a)=\pi(b) \in X^{(n)}_{\nu}
\} } \in Z_{2n} ( X^{[n]} \times X^{[n]} ).
$$
Note that the image of $D^{\nu}$ under either projection is the
closed stratum $\overline X^{[n]}_{\nu}$.
\begin{rmk}
\label{componts} 
{\rm $D^{\nu}$ is an irreducible component of  
$(\overline X^{[n]}_{\nu} \times_{X^{(n)}} \overline X^{[n]}_{\nu})_{red}$.
The latter 
contains  all the $D^{\mu}$ such that $\mu \preceq \nu$
precisely as its
 irreducible components.} 
\end{rmk}

Let $\nu \in \prt$.
The components of $p_{12}^{-1}(\, \!^t\hat\Gamma^{\nu}) \cap
p_{23}^{-1}( \hat\Gamma^{\nu} )$ in $X^{[n]} \times X^{(\nu)} \times X^{[n]}$
are all  of the expected dimension. 
The refined formalism mentioned above, 
 Remark \ref{propint} and Remark \ref{componts} show that  the following
is a well-defined $2n$-dimensional cycle with rational coefficients
$$
\Delta_{\nu} : = \frac{1}{m_{\nu}} \; \hat\Gamma^{\nu} \circ \!^t\hat\Gamma^{\nu}
 \quad \in Z_{2n} ( X^{[n]} \times X^{[n]} ) \otimes_{\zed} \rat
$$
 with the property that 
$\Delta_{\nu}$ is supported precisely on 
$(\overline X^{[n]}_{\nu} \times_{X^{(n)}} \overline X^{[n]}_{\nu})_{red}$
and that
$$
\Delta_{\nu} = \sum_{\nu \succeq  \nu'}{\epsilon_{\nu'}^{\nu} D^{\nu'}},
$$ 
where  the numbers $\epsilon_{\nu'}^{\nu}$ are   nonzero rational numbers
with the same sign as $m_{\nu}$.

\begin{lm}
\label{projd}
Let $\nu$ and $\mu$ be partitions of $n$. Then
$$
\Delta_{\nu} \circ \Delta_{\mu} = \delta_{\nu \mu} \Delta_{\nu}
\quad \in Z_{2n} ( X^{[n]} \times X^{[n]} ) \otimes_{\zed} \rat.
$$
where $\delta_{\nu \mu}$ is the usual Kr\"onecker function.

\noindent
In particular, the $\Delta_{\nu}$ are mutually orthogonal projectors
and:
$$
{ \Delta_{\nu} }_*: A_{\ast} ( X^{[n]} ) \to
A_{\ast} ( X^{[n]} ), \qquad  \qquad { \Delta_{\nu} }_* \circ { \Delta_{\mu} }_*=
\delta_{\mu \nu} { \Delta_{\nu} }_*.
$$
\end{lm}
{\em Proof.}
By the associativity of the composition of correspondences  we have
$$
\Delta_{\nu} \circ \Delta_{\mu} = \frac{
\hat{ \Gamma}^{\nu} \circ    \!^t\hat{\Gamma}^{\nu}  }{m_{\nu}} 
\circ
\frac{ \hat{ \Gamma}^{\mu} \circ    \!^t\hat{\Gamma}^{\mu}  }{m_{\mu}}=
\frac{ \hat{ \Gamma}^{\nu} }{ m_{\nu} } 
\circ \left(
 \frac{ \!^t\hat{\Gamma}^{\nu} \circ  \hat{ \Gamma }^{ \mu} }{ m_{\mu} }  \right)
 \circ \!^t\hat{\Gamma}^{\mu}.
$$
We conclude by Proposition \ref{ab} and Proposition \ref{aa}. \blacksquare

\bigskip

Recall that  $F_{\mu}$ denotes the  reduced  fiber
$(\pi^{-1}(\underline{x}))_{red}$, where $\underline{x} \in X^{(n)}_{\mu}$.
\begin{lm}
\label{df}
Let $\nu$ and $\mu$ be partitions of $n$.

(i) If $\nu \nnn \,\mu$, then ${\Delta_{\nu}}_* ([F_{\mu}])=0$. If
  $\nu =\mu$, then  ${\Delta_{\nu}}_* ([F_{\mu}])= $  $[F_{\mu}]$. 

(ii) If $\nu \nnn  \, \mu$, then   ${D^{\nu}}_* ( [ F_{\mu} ] )=0$. 
If $\nu =\mu$, then  ${\Delta_{\nu}}_* ([F_{\mu}])=$ 
$c_{\nu} [F_{\nu}]$,  where $\rat \ni c_{\nu} \neq 0$.
\end{lm}
{\em Proof.} 
We  compute ${\Delta_{\nu}}_*$
using refined intersections so that the classes $[F_{\mu}]$, which are non zero in
$A_{\ast}(F_{\mu})$, may be zero in $A_{\ast}(X^{[n]})$. Of course this does not happen if, for example, $X$ is projective.

  If $\nu$ does not refine $\mu$, then
$p^{-1}(F_{\mu}) \cap \Delta_{\nu} =\emptyset$, and the first parts of (i)
and (ii)  follow at 
once.

Note that since we have proved the first parts of (i) and (ii),
(i) implies (ii) by consideration of supports. 
We now prove  the second part of (i).
Let $\mu=\nu$. 
Note that ${ \Delta_{\nu} }_* ([F_{\nu}])=$ $\frac{1}{m_{\nu}} \hat\Gamma^{\nu}_*( \,\!^t \hat\Gamma^{\nu}_*(
[F_{\nu}]))$. The zero cycle $\; \!^t \hat\Gamma^{\nu}_*(
[F_{\nu}])$ is supported at the  point $\nu^{-1}(\underline{x}) \in X^{(\nu)}$. 
By virtue of the projection formula, its degree
equals $\deg{ ( [ F_{\nu} ] \cdot 
[  \overline X^{[n]}_{\nu}   ] ) } = m_{\nu}$; see Lemma \ref{ellings}.
In other words
$\!^t \hat\Gamma^{\nu}_*(
[F_{\nu}]) = m_{\nu}[\nu^{-1}(\underline{x}) ]$. We apply $\hat\Gamma^{\nu}_*$ and we find the conclusion. 
\blacksquare

\bigskip
Note that  the diagonal $\Delta=D^{1^n}$. 
The main reason for   introducing the correspondences
$\Delta_{\nu}$ is the following

\begin{pr}
\label{identity} The map
$
(\sum_{\nu \in \prt}{\Delta_{\nu}} - \Delta)_*
$
is the zero map.
In particular,
$
\sum_{\nu}{\Delta_{\nu}}_* = \Delta_* = Id_{A_{\ast}(X^{[n]})}.
$
\end{pr}
{\em Proof.}
This follows immediately from Theorem \ref{maintm}
which implies that $\!^t\hat{\Gamma}_*$ is the inverse of $\hat{\Gamma}_*$
and from the fact that $  \hat{\Gamma}_* \circ   \!^t\hat{\Gamma}_*=
\sum_{\nu \in \prt}{{\Delta_{\nu}}_*}$.
\blacksquare
\bigskip

\begin{pr}
\label{disd}
$
\sum_{\nu}{\Delta_{\nu}} - \Delta =0$ in   $ Z_{2n}(X^{[n]} \times X^{[n]})\otimes_{\zed}
\rat.
$
\end{pr}
{\em Proof.}
Without loss of generality, we may assume that $X$ is projective.
Note that, for every $\nu \neq 1^n$,  $D^{\nu}_*$ is identically zero
on $A_0(X^{[n]})$.  It follows,
by virtue of Lemma \ref{df},  that we
can   write $\sum_{\nu \in \prt}{\Delta_{\nu}} - \Delta= \sum_{\nu \neq 1^n}{p_{\nu}D^{\nu}}$,
where $p_{\nu}$ are suitable rational numbers.
Seeking a contradiction, assume that there is at least one partition $\mu$
for which $p_{\mu} \neq 0$. Choose a partition $\theta$ 
which is maximal, with respect to the partial
ordering on partitions, among all partitions $\mu$ for which $p_{\mu} \neq 0$.
We have, again by virtue of  Lemma \ref{df}, that
$$
0= (\sum_{\nu}{p_{\nu}D^{\nu}})_*[F_{\theta}]=p_{\theta}D^{\theta}_* (
[F_{\theta}] ) \neq 0.
$$
The contradiction stems from the projectivity assumption, for then no effective cycle 
can be trivial.
\blacksquare

\subsection{The structure and the generating function of  motives}
Let $X$ be an irreducible nonsingular projective surface defined over an algebraically closed field.
The correspondence $\hat{\Gamma}^{\nu}$ defines a morphism, which by
abuse of notation we denote by the same symbol,
of effective Chow motives with rational coefficients:
$$
\hat{\Gamma}^{\nu}: (X^{(\nu)}, \Delta_{X^{(\nu)}}) (n-l(\nu)) 
\longrightarrow
(X^{[n]}, \Delta_{X^{[n]}}).
$$
See \ci{manin},  page 459. The Tate-type shift is in accordance with
the usual convention $(\pn{1}, \Delta_{\pn{1}})(1)$ $=(\pn{1}, p)$, where
$p= \pn{1}\times \{0\}$.
This morphism admits,  in the category of
Chow motives with rational coefficients, a right inverse. This is  
given, again by abuse of notation,
by $\, \!^t\hat{\Gamma}^{\nu}$. By virtue of \ci{manin}, page 453, we can split-off a direct summand:
$$
\hat{\Gamma}^{\nu}: (X^{(\nu)}, \Delta_{X^{(\nu)}}) (n-l(\nu)) \simeq (X^{[n]}, \Delta_{\nu}).
$$
By virtue of  Lemma \ref{projd}, the projectors $\Delta_{\nu}$ are mutually orthogonal so that  we can split-off 
 a direct summand for each partition:
$$
\hat{\Gamma}:= \bigoplus_{\nu  \in \prt}
\hat{\Gamma}^{\nu}: 
\bigoplus_{\nu \in \prt}(X^{(\nu)}, \Delta_{X^{(\nu)}}) (n-l(\nu)) \simeq 
\bigoplus_{\nu \in \prt}(X^{[n]}, \Delta_{\nu})=(X^{[n]}, \sum_{\nu \in \prt}{\Delta_{\nu}}).
$$
The following is the main theorem of this section and follows immediately
from what above and Proposition \ref{disd}.
\begin{tm}
\label{motiff}
Let $X$ be an irreducible nonsingular projective surface defined over an algebraically
closed field.
There is a natural  isomorphism of effective Chow motives with rational coefficients
$$
(X^{[n]}, \Delta_{X^{[n]}}) \simeq \bigoplus_{\nu \in \prt}(X^{(\nu)}, \Delta_{X^{(\nu)}}) (n-l(\nu)).
$$
\end{tm}
\begin{rmk}
\label{anyff}
{\rm 
Theorem \ref{motiff} and Theorem \ref{vmot} below hold over any field. We leave the necessary but easy  modifications to the reader.}
\end{rmk}
\begin{rmk}
\label{coho}
{\rm Over the complex numbers, 
 Theorem \ref{motiff} gives immediately the structure
of the singular cohomology $H^*(X^{[n]}, \rat)$, together with its Hodge structure.
The resulting isomorphisms coincide with the ones obtained in
 \ci{decatmig}.  
}
\end{rmk}
\bigskip
Voevodsky \ci{vo} has defined several motivic categories over a field
$k$: $DM^{eff}_{gm}(k)$
(cf. \ci{vo}, 2.1.1), $DM^{eff}_{-}(k)$ (cf. \ci{vo}, 3.1.12), $DM^{eff}_{-,et}(k)$ (cf. \ci{vo}, 3.3), $DM_{h}(k)_{\rat}$ (cf. \ci{vo} and \ci{vo1}).
These are, essentially, categories of bounded complexes of pre-sheaves on 
certain sites. There is an obvious formalism of shifts of complexes and a formalism of Tate-type shifts.
 In addition, there are natural transformations between
 the category of effective Chow motives
with rational coefficients and these other  categories (see \ci{vo}, 2.1.4):
$$
Chow^{eff}(k)_{\rat} \to 
DM^{eff}_{gm}(k) \to 
 DM^{eff}_{-}(k) \to  DM^{eff}_{-,et}(k)  \to 
 DM_{h}(k)_{\rat}.
$$
The formalism holds for quotient varieties.

Let us denote by $M(Z)$ the object in any of the above categories corresponding to
a variety $Z$.
Theorem \ref{motiff} and the natural transformations mentioned above imply
the following
\begin{tm}
\label{vmot}
Let $X$ be an irreducible nonsingular projective surface defined over an algebraically closed field. Then
$$
M(X^{[n]}) \simeq \bigoplus_{\nu \in \prt} M(X^{(\nu)}) (n-l(\nu))[2n -2l(\nu)].
$$
\end{tm}
If we use the dual Tate twist, then we can re-write the conclusion of the theorem as follows:
$$
M(X^{[n]})(n)[2n] \simeq \bigoplus_{\nu \in \prt} M(X^{(\nu)}) (l(\nu))[2l(\nu)].
$$
This is formally reminiscent of 1) the decomposition theorem for the Douady-Barlet morphism
\ci{decatmig}, where shifts for complexes of sheaves appear and 2) the Hodge structure of
$X^{[n]}$ (see \ci{decatmig}, for example), where Tate shifts appear. It is by staring at those
two formulas that we became convinced that the statement of Theorem 
\ref{vmot} should hold. Theorem \ref{maintm}  and Theorem 
\ref{motiff} were the means to realize this isomorphism.

\bigskip
It is amusing to note that there is a generating function for this picture, the one
that generates partitions. In fact, skipping the bookeeping details, the formula is as follows:

$$
\sum_{n=0}^{\infty}{X^{[n]}(n)[2n]} \, q^n
= \prod_{m=1}^{\infty}{ \frac{1}{ ( 1-X_m(1)[2]\, q^m  ) } },
$$
where it is understood that in the r.h.s. we  identify
monomials of type $X_1^{a_1} \cdots X_t^{a_t}$ with
$X^{(\nu)}$, where $\nu = 1^{a_1} \cdots t^{a_t} \in \frak{P}(t)$.
Note that this suggests a series of identities in the 
Grothendieck ring of varieties, among cohomology groups, mixed Hodge structures,
intersection cohomology, equivariant $K$-theories, Chow groups, motives,
 etc., all of which are true, except, possibly, for the first one.

Authors's addresees:

\smallskip
Mark Andrea A. de Cataldo,
Department of Mathematics,
SUNY at Stony Brook,
Stony Brook,  NY 11794, USA. \quad
e-mail: {\em mde@math.sunysb.edu}

\smallskip
Luca Migliorini,
Dipartimento di Matematica , Universit\`a di Trento,
Via Sommarive, 14,
38050 Povo (Tn),  ITALY. \quad
e-mail: {\em luca@alpha.science.unitn.it}


\begin{thebibliography}{99}




\bibitem{br}{J. Brian\c{c}on, ``Description de $Hilb^nC\{x,y\}$," Invent. Math.
{\bf 41} (1977), 45-89.
}

\bibitem{ch}{J. Cheah, ``The Hodge numbers  of the Hilbert scheme of points of
 a smooth projective surface," J. Algebraic Geometry {\bf 5} (1996), 479-511.}


\bibitem{decatmig} M. de Cataldo, L. Migliorini, ``The Douady space
of a complex surface," to appear in Adv. in Math.


\bibitem{e-s1}{G. Ellingsrud, S.A. Stromme,  ``On a cell decomposition of the 
Hilbert scheme
of points in the plane," Invent. Math.  {\bf 91} (1988) 365-370.}

\bibitem{ellstrom} {G. Ellingsrud, S.A. Stromme,  ``An intersection number
for the punctual Hilbert scheme of a surface," Trans. Amer. Math. Soc.
{\bf 350} (1998), no.6, 2547-2552.}

\bibitem{fo}{J. Fogarty, ``Algebraic families on an algebraic surface," 
American Journal
of Math. {\bf 90} (1968), 511-521.}

\bibitem{fulton} W. Fulton, {\em Intersection Theory},
 Ergebnisse der Mathematik,
3.folge. Band 2, Springer-Verlag, Berlin Heidelberg 1984.



\bibitem{go}{L. G\"ottsche, {\em Hilbert schemes of zero-dimensional subschemes
of smooth varieties}, LNM 1572, Springer-Verlag Berlin 1994.
}


\bibitem{go-so}{L. G\"ottsche, W. Soergel, ``Perverse sheaves and the cohomology of the
Hilbert schemes of smooth algebraic surfaces," Math. Ann. {\bf 296} (1993), 235-245.}




\bibitem{gr}{L. Grojnowski, ``Instantons and affine algebras, I. The Hilbert scheme and 
vertex operators," Math. Res. Letters {\bf 3} (1996), no. 2, 275-291}


\bibitem{ia}{A. Iarrobino, {\em Punctual Hilbert schemes}, Mem. Amer. Math.
Soc. {\bf 188}, 1977.}

%\bibitem{le}{M. Lehn, ``Chern classes of tautological sheaves on Hilbert schemes 
%of points on surfaces," 
%Invent. Math. 136 (1999), no. 1, 157--207.}



\bibitem{manin} Ju. I. Manin, 
``Correspondences, motifs and monoidal transformations," 
Math. USSR-Sb. 6 (1968), 439--470.




\bibitem{nakp}{H. Nakajima, ``Heisenberg algebra and Hilbert schemes of points
on projective surfaces," Ann. of Math. (2) {\bf 145} (1997), no.2, 379-388.}


\bibitem{nak}{H. Nakajima, {\em Lectures on Hilbert schemes of points on
surfaces}.  University Lecture Series, 18. American Mathematical Society,
Providence, RI, 1999.}





\bibitem{va-wi}{Vafa, E. Witten, ``A strong coupling test for S-duality,"
Nucl. Phys. {\bf 431} (1994), 3-77.}

\bibitem{vi}{A. Vistoli, ``Higher equivariant $K$-theory for finite group actions,"
Duke Math. J. {\bf 63} (1991), no.2, 399--419.}


\bibitem{vo} V. Voevodsky, ``Triangulated categories of motives over a field," 
preprint, k-theory arch. n.74.

\bibitem{vo1} V. Voevodsky, ``Homology of schemes, I," preprint, k-theory arch. n.31.


\end{thebibliography}
\end{document}